\documentstyle{article}
\newcommand{\Tilde}{\widetilde}
\newcommand{\R}{{\bf R}}

\newcommand{\Z}{{\bf Z}}

\newcommand{\al}{{\alpha}}

\newcommand{\be}{{\beta}}

\newcommand{\om}{{\omega}}
\newcommand{\eps}{{\varepsilon}}

\newcommand{\ga}{{\gamma}}
\newcommand{\Ga}{{\Gamma}}
\newcommand{\ka}{{\kappa}}
\newcommand{\la}{{\lambda}}
\newcommand{\si}{{\sigma}}
\newcommand{\Symp}{{\rm Symp}}
\newcommand{\Ham}{{\rm Ham}}

\newcommand{\Xx}{{\cal X}}

\newcommand{\Q}{{\bf Q}}

\newcommand{\Si}{{\Sigma}}

\newcommand{\1}{{\rm id }}

\newcommand{\SS}{{\smallskip}}
\newcommand{\MS}{{\medskip}}
\newcommand{\BS}{{\bigskip}}
\newcommand{\NI}{{\noindent}}
\newcommand{\proof}[1]{\noindent{\bf Proof#1:\  }}

\newcommand{\QED}{\hfill$\Box$\medskip}
\newtheorem{theorem}{Theorem}[section]
\newtheorem{cor}[theorem]{Corollary}

\newtheorem{remark}[theorem]{Remark}
\newtheorem{lemma}[theorem]{Lemma}

\newtheorem{prop}[theorem]{Proposition}

\begin{document}

\title{On the Flux Conjectures}
\author{Fran\c{c}ois Lalonde\thanks{Partially supported by NSERC grant 
OGP 0092913
and FCAR grant ER-1199.} \\ Universit\'e du Qu\'ebec \`a Montr\'eal
\\ (flalonde@math.uqam.ca) \and Dusa McDuff\thanks{Partially
supported by NSF grant DMS 9401443.} \\ State University of New York
at Stony Brook \\ (dusa@math.sunysb.edu)
 \and Leonid
Polterovich\thanks{All three authors supported by a NSERC
Collaborative Project Grant CPG0163730.} \\
Tel-Aviv University \\ (polterov@math.tau.ac.il)}

\date{February 18, 1997}
\maketitle

\begin{center}
{\bf Abstract}
\end{center}

The ``Flux conjecture'' for symplectic manifolds states that the group of
Hamiltonian  diffeomorphisms is $C^1$-closed in the group of all
symplectic diffeomorphisms.
We prove the conjecture for {\it spherically rational} manifolds and for
those whose minimal Chern number on 2-spheres either vanishes or is 
large enough. We also confirm  a natural version of the Flux conjecture
for symplectic torus actions.
In some cases we can go further and prove that the group of Hamiltonian 
diffeomorphisms is $C^0$-closed in the identity component of the
group of all symplectic diffeomorphisms.

\section{Introduction}

   Let $(M, \om)$ be a compact symplectic manifold without boundary, 
and $G = \Symp_0(M,\omega)$
be the identity component 
of the group of symplectic diffeomorphisms of $M$.
There is an exact sequence
$$
\Tilde H \to \Tilde {G}
\stackrel{F}{\to} H^1(M; \R) $$
where $\Tilde G$ is 
the universal cover of $G$, $\Tilde H$ is the subgroup of $\Tilde G$
formed by those paths which are homotopic in $G$ with fixed endpoints
to a Hamiltonian path, and $F$ is the flux homomorphism defined by
$$
F(\{\phi_t\}) = \int_0^1 \la_t \, dt.
$$
Here the  family of closed $1$-forms $\{\la_t\}$  generates
 the isotopy $\{\phi_{t\in [0,1]}\}$ in the sense that
$$
\om(\frac{d\phi_t}{dt}, \cdot) = \la_t(\cdot)\quad\mbox { for all }\;\;t\in
[0,1]. 
$$
  Recall 
also that a path
is Hamiltonian if its generating family $\{\la_t\}$  is exact
at each time $t$.
 Thus the exact sequence simply
expresses the fact that a path with vanishing (average) flux can be 
perturbed to have vanishing ``instantaneous
flux'' for each $t$. 
Moreover, it is easy to check that $F(\{\phi_{t}\})$ is the class in 
$H^1(M;\R)$ that
assigns to each loop $\ga$ in $M$ the integral of $\om$ over the cylinder
$$
C_\ga : S^1 \times [0,1] \to M 
$$
defined by $C_\ga(z,t) = \phi_t(\ga(z))$.
(Proofs of the above statements
can be found in Banyaga~\cite{B} or
McDuff--Salamon \cite{MS}.)  We will sometimes refer to the cylinder
$C_\ga$ as the {\em trace} of $\{\phi_t\}$ on $\ga$.  

   Let us denote by $\Ham (M,\om)$, or simply $\Ham(M)$, the group of
Hamiltonian diffeomorphisms, that is to say those elements of $G$ which
are endpoints of Hamiltonian paths.\footnote
{
When no specific mention is made, all our paths begin
at the identity.  Note also that the group $\Tilde H$
which appears in the above exact sequence is just the universal cover of
$\Ham(M,\om)$.} 
Further, we define the {\it Flux subgroup} $\Ga$ of $H^1(M; \R)$
to be the image by $F$ of the closed loops in $G$: 
$$
 \Ga = {\rm Image}\left( \pi_1(G) \stackrel{F}{\to} H^1(M; \R) \right).
$$
If $ \{\phi_{t \in [0,1]}\}$ is a path of symplectic
diffeomorphisms then it is easy to see that  its endpoint
$\phi_1$ belongs to $\Ham (M)$ if and only if its flux $F(\{\phi_t\})$
belongs to $\Ga$.
In other words, there is an exact sequence
$$
(*)\qquad\qquad\Ham(M,\om)\to \Symp_0(M,\om)\to
H^1(M,\R)/\Ga, $$
where the second map is induced by $F$.
  Thus, the
 group $\Ga$  contains
crucial information about the manifold, being the key
to whether or not a symplectic diffeomorphism is Hamiltonian. 
In particular it is important
to know if $\Ga$  is discrete, and if  it is, how to estimate the
size of a neighborhood of $\{0\}$ in $ H^1(M;\R)$ that contains no element
of $\Ga$ except $\{0\}$ itself. We begin with the following observation:

\begin{prop}\label{prop:equivalence} For any closed symplectic manifold,
$\Ga$ is discrete if and only if  the subgroup of Hamiltonian
diffeomorphisms is $C^1$-closed in the full group of (symplectic)
diffeomorphisms of the manifold. \end{prop}

This follows almost immediately from the existence of the exact
sequence $(*)$ since all the maps involved are $C^1$-continuous. 
A more detailed proof is given in \S3 below. 
   We can now state the main conjectures:
\MS

\NI
{\bf The $C^1$-Flux Conjecture.} 
{\it The flux subgroup is discrete for all 
symplectic manifolds}. 

\MS

\NI
{\bf The $C^0$-Flux Conjecture.} {\it For any symplectic manifold, the group 
$\Ham(M)$ is $C^0$-closed in the identity component of the group of 
symplectic diffeomorphisms.}

\BS

 The main  
problem posed by the Flux conjecture\footnote{Classically, 
the ``Flux conjecture'' is what
we call here the $C^1$-Flux conjecture.  For this reason, when we refer to the
``Flux conjecture'' with no qualification, we always have the $C^1$-case in
mind.
}
   can be presented in the following
way.  Consider a ``long'' Hamiltonian
path $\{\phi_{t \in [0,1]}\}$ beginning at the identity whose graph at 
intermediate times escapes out of a Weinstein 
neighborhood $U$ of the diagonal in $(M \times M, -\om \oplus 
\om)$, but which comes  back inside that neighborhood at time $t=1$. 
As is well-known, $U$ can be identified with a neighborhood of the zero
section in the cotangent bundle of $M$, and, when it is,  a Hamiltonian
isotopy that stays inside $U$ is given by the
graphs of a family of exact $1$-forms.  However, if the path goes outside
$U$, then in
general there is no obvious reason that the $1$-form corresponding to
its endpoint be   exact.

\MS

Note also that the statement that $\Ham(M)$ is $C^k$-closed for one $k
\ge 1$ is obviously equivalent to the statement that it is $C^k$-closed for
all $k\ge 1$.  However, the question of whether it is $C^0$-closed is a
priori quite different, and  seems to be very difficult to decide except in some
particular cases like tori (see below).  This is a fundamental issue.
Indeed all known invariants (such as symplectic homology) are
$C^0$-invariants attached  to Hamiltonian {\it paths}, but one is actually
interested in the dependence on the  Hamiltonian {\it endpoint}. 
Observe also that, although symplectic rigidity tells us that the
group $\Symp(M)$ of all symplectomorphisms is $C^0$-closed in the group of
all diffeomorphisms of $M$, it is not known whether  or
 not the identity component
$\Symp_0(M)$ of $\Symp(M)$ is $C^0$-closed in $\Symp(M)$.
To avoid this question, we look here at the $C^0$-closure of
$\Ham(M)$ in $\Symp_0(M)$.

As 
will become clear in this paper, the Flux conjectures lie at the very
borderline between soft and hard symplectic topology. What we do 
here is to show how hard techniques and their consequences
can be used to prove these 
conjectures in many cases where purely soft, topological
methods seem to fail.

\subsection{Some old results}\label{old}

   The flux  group $\Ga$ was first explicitly mentioned
 in Banyaga's foundational paper~\cite{B}, where it was observed
that $\Ga$ is discrete when $M$ is K\"ahler, or, more generally, when
$M$ is Lefschetz.\footnote{
This means that 
$
H^1(M,\R) \stackrel{\cup [\om]^{n-1}}
{\longrightarrow} H^{2n-1} (M;\R)
$
is an isomorphism.}  For the sake of completeness,
we will briefly recall this and other early results on $\Ga$.
\MS

\NI
{\bf (i)} It is easy to check that 
there is a commutative diagram expressing the fact
that the following
two homomorphisms coincide up to a universal
multiplicative constant depending only on the dimension of $M$:
$$
\pi_1(G) \stackrel{ev}{\to} \pi_1(M) \to H_1(M;\Z) \stackrel{PD}{\to}
   H^{2n-1} (M;\Z)
$$
and
$$
\pi_1(G) \stackrel{F}{\to} H^1(M,\R) \stackrel{\cup [\om]^{n-1}}{\to} 
     H^{2n-1} (M;\R).
$$
Here $ev$ is the evaluation map at a point of $M$, $PD$ is the 
Poincar\'e dual and $F$ is the flux homomorphism.   
But the first map has discrete image. Thus, if a manifold is Lefschetz,
the subgroup $\Ga \subset H^1(M,\R)$ is discrete too.  Note also that
in this case the flux homomorphism vanishes 
on all elements of $\pi_1(G)$ that
evaluate trivially in the homology group $H_1(M,\R)$.\MS

\NI
{\bf (ii)} Temporarily, let $G$ denote the group of diffeomorphisms 
which preserve some
given structure on $M$. Then the ``$c$-flux'' homomorphism 
$$
\pi_1(G) \stackrel{F_c}{\to} H^1(M,\R),
$$ 
defined by assigning to the pair  
$(\phi_{t \in [0,1]},\ga) \in \pi_1(G) \times H_1(M;\R)$ 
the integration of the class $c \in H^2(M)$ 
on the trace of $\{\phi_t\}$ on $\ga$, 
vanishes identically if $c$ is a characteristic 
class for the structure preserved by $G$.  This was proved by McDuff 
in~\cite{McDF} by an argument based on
Gottlieb theory.  (A more transparent proof will be given
in~\cite{LMP}.)  Restricting  to the case when $G$ is the group of
symplectic diffeomorphisms and $c$ is the first Chern class of the
tangent bundle of the symplectic manifold, we see that the flux vanishes 
for  monotone
manifolds.\footnote{
We recall that a {\it monotone}
manifold is a symplectic manifold satisfying $c = \ka \om$ for some
$\ka > 0$.  It is {\it semi-monotone} if $\ka \ge 0$.}
(This
application to monotone manifolds was not explicitly  mentioned
in~\cite{McDF}, but appears in Lupton--Oprea~\cite{LO}.)  \MS

\NI
{\bf (iii)} It is also proved in~\cite{McDF}
that $\Ga$ is discrete
 when the class of the symplectic form is 
decomposable, that is to say when $[\om]$
is the sum of products of elements of $H^1(M,\R)$.
\MS

\NI
{\bf (iv)}  Finally it follows immediately from
Ginzburg's results in~\cite{GI} that any symplectic  
torus action on a Lefschetz manifold has discrete flux: see 
\S\ref{se:autonomous} for more detail.

\subsection{The evaluation map}

In this section we discuss the relation of the Flux conjectures to the 
topological properties of the orbits $\{\phi_t(x)\}$ of a loop
$\{\phi_t\}$.  Our work is based on the following
deep fact:
\begin{quote}\em For every Hamiltonian
flow $\{\phi_{t\in [0,1]}\}$ on a closed symplectic manifold, there is at
least one fixed point $x\in M$ of $\phi_1$  such that the
loop $\{\phi_t(x)\}_{t\in [0,1]}$ is   {\it contractible}. In particular,
any  loop in $G$ with vanishing flux is sent by the evaluation map
 to $0 \in
\pi_1(M)$.
\end{quote}
This is proved along with  the Arnold conjecture:  indeed
all proofs of the Arnold conjecture  give a nonzero lower bound for
the number of fixed points of $\phi_1$ with {\em contractible} orbits
$\{\phi_t(x)\}$: see
Hofer--Salamon~\cite{HS},
Fukaya--Ono~\cite{FO} and Liu--Tian~\cite{LIUT}.  

As we will see in \S3, the following result is an almost immediate consequence.

\begin{prop} \label{prop:commutative} 
The  flux homomorphism has discrete
image on $\pi_1(G)$  if and only if it has discrete image on the
subgroup $K$ of $\pi_1(G)$ formed by loops with contractible orbits in
$M$.     \end{prop}

Thus the Flux conjecture only depends on those symplectic
loops which, like Hamiltonian ones, have contractible orbits. 
One sometimes says that the flux homomorphism of some symplectic
manifold {\it factorizes} if it factorizes through the
evaluation map $\pi_1({\rm Symp}(M)) \to \pi_1(M)$.  In other
words, it factorizes when every symplectic loop with contractible
orbit has  vanishing flux. Proposition~\ref{prop:commutative}
shows that factorization implies discreteness. 
We will discuss a result of Th\'eret's 
on  factorization at the end of \S1.4.

\begin{remark}\rm  Although both the above statements remain true for
compactly supported Hamiltonian flows on noncompact manifolds,
they do not extend to flows of
arbitrary support, even if one considers the evaluation in {\it homology}. Indeed,
the standard rotation of the annulus is Hamiltonian. On the other hand the fact
that the evaluation in real homology vanishes for Hamiltonian loops on closed
manifolds is very elementary and follows from the commutative diagram in \S1.1.
This illustrates a striking contrast between the homological (mod torsion) and
homotopical points of view on the evaluation map:
see Bialy-Polterovich~\cite{BPGAFA} for further discussion.\QED
\end{remark}

In view of the above results, one might wonder if any symplectic loop
with contractible orbits has to be Hamiltonian up to homotopy.  
In dimension $4$ the
answer to this question is not known in general, but it is false in general.
Indeed, one has the following result:  

 \begin{prop}\label{prop:orb}  Let $(M,\om)$ be a closed symplectic manifold.

\NI
(i) In dimension $4$ every symplectic $S^1$-action
with contractible orbits is Hamiltonian (and hence  has fixed points.)

\NI
(ii)  There is a nonHamiltonian symplectic $S^1$-action on a $4$-manifold
with orbits that are homologous to zero but not
contractible.

\NI
(iii)  In dimension $6$ there is a nonHamiltonian $S^1$-action
with  fixed points and hence with contractible orbits.
\end{prop}

We will see in \S2 that part (i) follows almost immediately from results in
McDuff~\cite{McD:mom}.  It is easy to construct an example 
for (ii) in which  $M$ 
is a nonK\"ahler $T^2$-bundle over
$T^2$,  for example  the  Kodaira--Thurston
 manifold, while an
example of type (iii)  may be found
in~\cite{McD:mom}.   

The above result shows that
 it is not in general possible to distinguish Hamiltonian
loops from others by looking only at the topological 
properties of their orbits.  
However, the story becomes more interesting if one looks at the
 ``evaluation map'' $ev_\phi: H_k(M) \to H_{k+1}(M)$
of a  loop $ \{\phi_{t \in [0,1]}\}$
in dimensions $k > 0$.  Here, $ev_\phi$ is the map which
takes a $k$-cycle $\ga$ to its trace  $C_\ga =
\cup_t\{\phi_t(\ga)\}$.  Thus  $\{\phi_t\}$ is Hamiltonian
up to homotopy
if and only if $[\om]$ vanishes on the image of $ev_\phi|_{H_1}$.

Recently we have shown that when $M$ is monotone $ev_\phi$
vanishes identically for all $k$.  The
proof is based on Seidel's description in~\cite{SEI} of the canonical
action of $\pi_1(\Ham(M))$ on the quantum cohomology ring of $M$.
The vanishing of this map for $k=1$ 
means that the flux for Hamiltonian loops
$$
H_1(M) \stackrel{ev_\phi}{\to} H_2(M) 
                          \stackrel{\int \om}{\to}  \R
$$  
actually vanishes for {\it topological reasons}. Thus, if a loop
$\be= \{\phi_{t}\}$  in ${\rm Diff}(M)$ is Hamiltonian with respect to $\om$
and is homotopic to a loop $\be'$ which is symplectic with respect to
some other symplectic form $\om'$, then $\be'$ is necessarily {\it
Hamiltonian} (up to homotopy) with respect to $\om'$. This implies that if a 
loop is
Hamiltonian with respect to $\om$, then, for any symplectic form $\om'$
sufficiently close to $\om$, it can be homotoped to a $\om'$-symplectic loop
which is actually Hamiltonian. Note that this result can be interpreted as an
obstruction in the following way. Let $\be$ be
any $\om$-symplectic nonHamiltonian loop  and 
denote by $A_{\be,\om} \neq \{0\}$
the image of the map $H_1(M,\Z) \to H_2(M,\Z)$ induced by $\be$. Let $a \in
H^2(M,\R)$ be any class that vanishes on $A_{\be,\om}$. Then a deformation of
the symplectic  structure from the class $[\om]$ to $a$ -- if it exists -- 
cannot be lifted to a deformation of the loop $\be$. In other words, 
each space $A_{\be,\om}$ gives an
obstruction either to the deformation of $\om$ or to the deformation of the
image of $\pi_1({\rm Symp}(M,\om))$ inside $\pi_1({\rm Diff}(M))$. 
Because a symplectic structure $\om$ and an $\om$-symplectic loop $\be$
determine a symplectic fibration $V$ with fiber $(M,\om)$ over the
$2$\/-sphere, this can also be interpreted in the following way.  The
existence of a ruled symplectic form on $V$ compatible with a given
symplectic fibration $V\to S^2$ depends  only on the connected component
of that fibration in the space of all symplectic fibrations. These results and
various generalizations and corollaries will appear in our forthcoming
paper~\cite{LMP}.

\subsection{The $C^1$-Flux conjecture}

We now present a list of manifolds for which we have been
 able to confirm the $C^1$-Flux conjecture.

We start our discussion with the
the case of symplectic torus actions. To motivate it,
suppose that, for some symplectic manifold $M$, the group $G$ retracts
onto a finite dimensional Lie subgroup $H$.
(This is known to be true for some``simple" symplectic manifolds,
for instance for surfaces and some of their products.)
Then all elements  of the fundamental 
group of $G$ are represented by elements of the fundamental
group of a maximal torus inside $H$, and the flux conjecture 
reduces to the same conjecture about the fluxes of an autonomous
action of the torus. Our first result confirms the Flux Conjecture
in this case. 

\begin{theorem} \label{thm:autonomous} Let $T^n$ act symplectically on a 
closed symplectic manifold. Then the restriction of the flux
homomorphism $F$ to $\pi_1(T^n)$ has discrete image in
$H^1(M,\R)$. \end{theorem}

The proof, which is given in \S2, combines ideas from Ginzburg~\cite{GI}
with  an analysis of the Morse-Bott singularities of the  corresponding
(generalized) moment map. 

Let us go back to the general non-autonomous case. Note that the 
conjecture obviously holds if the integration morphism
$\int \om: H_2^T(M, \Z) \to \R$ has discrete image, where
$H_2^T(M, \Z)$ denotes the set of classes that can be represented by
continuous maps of the $2$-torus. 
The first statement below
shows that it is enough to require that this
integration morphism has discrete image on {\it spherical classes}
$$
H^S_2(M, \Z) = {\rm Im}(\pi_2(M) \to H_2(M, \Z))
$$
alone. The second statement refers to the minimal spherical
Chern number, which is by definition the  nonnegative generator
of the image of $H^S_2(M, \Z)$ 
by the first Chern class $c = c_1(TM)$ of the tangent bundle of $M$.

\begin{theorem} \label{thm:main-C1} The $C^1$-Flux conjecture holds 
in  the following cases:
\begin{description}
\item[(i)] The manifold 
is spherically rational, that is to say the image
 of the spherical $2$-classes $H_2^S(M, \Z)$ 
by the integration morphism $\int \om$ is a discrete subgroup 
of $\R$ (of the form $\la \Z$ for some
nonnegative real number $\la$);  
\item[(ii)] The minimal spherical Chern number 
is either zero or is no less than $2n = {\rm dim}_{\R} (M)$; 
\item[(iii)] $M$ has dimension $4$ and $\pi_1(M)$ acts trivially on
$\pi_2(M)$. \end{description}
\end{theorem}

The various parts of this theorem have different proofs.  Parts (i) and (iii)
follow from
Proposition~\ref{prop:commutative} by easy topological 
arguments.   The main idea in the proof of part (ii) is the
following. Suppose that some class $a \in H^1(M, \R)$ belongs to
$\Ga$. Then 
the autonomous path $\{\psi_t\}$, generated by any representative $\la_a$
of $a$  has flux equal to $a$.
But $\psi_1$ is also the
endpoint of a Hamiltonian path $\{\phi_t\}$.  Thus one can compare
the Floer--Novikov homology of a generic perturbation of $\{\psi_t\}$ to
the Floer homology of the path $\{\phi_t\}$. Morally, this should lead to
strong constraints on the class $a$. In many cases, they are
strong enough to 
give a complete description of $\Ga$.  If one
is only interested in proving discreteness, one can limit this study
to the comparison 
of these two paths when $\la_a$ is so $C^1$-small  that
its zeros are  exactly the fixed 
points of $\psi_1$.
 In this case, one compares the ordinary Morse--Novikov 
homology of the small $1$\/-form $\la_a$ to some Floer homology. In particular,
if $\la_a$ 
does not have enough zeroes to satisfy the constraints
of the Arnold conjecture or if the zeroes do not have the appropriate 
indices, one concludes that it cannot belong to $\Ga$.

Theorem~\ref{thm:main-C1} extends the  results mentioned in \S1.1
above.  Indeed, an
 immediate corollary of (i) and (ii) is that the Flux conjecture
holds for semi-monotone manifolds.  
  Also (i) implies discreteness in the decomposable case, since any
decomposable form  vanishes on all 
spherical $2$-classes. 

\subsection{The $C^0$-Flux conjecture}

We complete the introduction with a 
discussion of the $C^0$-Flux conjecture.
Again, set $G=\Symp_0(M)$, and let $G_i \;(i=0,1)$ be
the $C^i$-closure of $\Ham(M)$ in $G$. 
Let $\Gamma_i$ be the image under the flux homomorphism of the
lift of $G_i$
to the universal cover of $G$. It is not hard to check that in this language
Proposition 1.1 states that ${\rm Closure}(\Gamma) = \Gamma_1$, while the
Flux conjecture and the $C^0$-conjecture are equivalent to the statements
$\Gamma=\Gamma_1$ and $\Gamma = \Gamma_0$ respectively.

It turns out that when $M$ is Lefschetz the group $\Gamma_0$
is contained in a group $\Gamma_{top}\subset H^1(M,\R)$
which depends only on the  topology of $M$ and on the
cohomology class of the symplectic form. Namely,
consider the space of all smooth (or even continuous) maps $M \to M$,
and denote by ${\rm Map}_0(M)$ the connected component of the identity.
The flux homomorphism $F$ extends naturally
to a homomorphism 
$\pi_1({\rm Map}_0(M))
\to H^1(M,\R)$ by integrating the
symplectic form on the evaluation of the path of maps
on the $1$\/-cycles of $M$. Denote by $\Ga_{top}$ its image
$F(\pi_1({\rm Map}_0(M)))$.  Clearly, 
$\Gamma \subset \Gamma_{0}$.

\begin{theorem}\label{thm:Lefschetz} If $M$ is Lefschetz, then
$
\Ga_0 \subset \Ga_{top} 
$.
\end{theorem}

See \S\ref{se:C0} for the proof of this
theorem.
Since $\Gamma \subset \Gamma_{0}$ we immediately get the following
useful consequence.  

\begin{cor} \label{thm:main-C0} Assume that $M$ is Lefschetz and
that $\Ga_{top} = \Ga$. 
Then the $C^0$-Flux conjecture holds for $M$.
\end{cor}  

As an example, take $M$ to be a closed K\"ahler manifold
of nonpositive curvature such that its fundamental
group has no center. For instance a product of surfaces
of genus greater than 1 has this property. It is easy to
see that in this case $\pi_1({\rm Map}_0(M))$ is trivial\footnote{
Note that  $\pi_1({\rm Map}_0(M))$ always maps into the {\it center} of
$\pi_1(M)$.}
and hence $\Gamma_{top} = \{0\}$. Moreover, $M$ is
Lefschetz since it is
K\"ahler. Thus the previous corollary implies that 
the $C^0$-Flux conjecture holds on $M$.

Another immediate application of this corollary is that
the $C^0$-Flux conjecture holds for
the $2n$-dimensional torus with a translation invariant
symplectic structure.
Indeed, in this case  $M$ is Lefschetz so that the Flux homomorphism
factors through $H_1(M,\Z)$.  Moreover
translations generate a large enough subgroup of
$\Symp_0(M)$ for us to see that $\pi_1(\Symp_0(M))$ maps onto $H_1(M,\Z)$.
Thus $\Ga = \Ga_{top}$, and so, by the above theorem,
the conjecture holds. We refer the reader
to McDuff--Salamon~\cite{MS} for the explicit computation
of $\Gamma$ for tori. 

 Another proof of the $C^0$-conjecture for these  tori
was first observed by Herman in 1983, just after Conley and Zehnder's
proof~\cite{CZ} of Arnold's conjecture for the  torus. His idea consists in the
following basic observation: if a path $\{\psi_t\}$ has flux $a$ and an endpoint
$\psi$ that is equal to the $C^0$-limit of Hamiltonian diffeomorphisms, then, by 
composing everything by an appropriate Hamiltonian,
we may assume that $\psi$ is the translation of flux
$a$.  We must show that $\psi$ is Hamiltonian.  But, since all Hamiltonian
diffeomorphisms  have a fixed point by Arnold's conjecture, so does the
limit $\psi$.  Hence  $\psi$ must be the identity, which, of course, is
Hamiltonian.

A natural generalization of Herman's idea is to replace
the counting of fixed points by an analysis of the
limiting behavior of Floer homologies associated
to a sequence of Hamiltonian diffeomorphisms, which is just the
approach we take in Theorem \ref{thm:main-C1} (ii).

Let us mention finally that the questions which we
discussed above can be posed in a more general
context of {\it Lagrangian submanifolds}
of a symplectic manifold. More precisely,
let $M$ be a symplectic manifold and $L \subset M$
be a closed Lagrangian submanifold. Denote by
$H(L)$ (resp. $S(L)$) the space of all Lagrangian submanifolds
which are obtained from $L$ by a Hamiltonian (resp. Lagrangian)
isotopy. The $C^i$-conjecture ($i=0,1$) in this case means
that $H(L)$ is $C^i$-closed in $S(L)$. A  result
in this direction was obtained recently by Theret~\cite{TH}.
See also Bialy--Polterovich~\cite{BPGAFA} for results on the
evaluation in $\pi_1(M,L)$ in this situation.
\MS

Throughout the paper we will assume that $M$ is a smooth 
compact manifold without boundary.  Many results apply when $M$ 
is noncompact, provided that we consider symplectomorphisms of
compact support.  We leave such extensions to the reader.

  We are grateful to M.~Braverman, M.~Faber, V.~Ginzburg, A.~Givental and 
J.-C.~Sikorav for
fruitful discussions on some aspects of this work.

\section{Torus actions} 
\label{se:autonomous}

 We begin by proving 
Theorem~\ref{thm:autonomous}, and then discuss the examples mentioned
in
Proposition~\ref{prop:orb}. 

\subsection{The discreteness of $\Ga$ for torus actions}

Throughout this 
section we use
the identification 
$$
\pi_1(T^d) = H_1(T^d,\Z) \subset H^1(T^d,\R).
$$

We start with the following useful notion introduced by
Viktor Ginzburg~\cite{GI}. A symplectic action of $T^d$ on $M$ is called
{\it cohomologically free} if the flux homomorphism
$F: H_1(T^d,\R) \to H^1(M,\R)$ is injective.
Clearly, for such an action  
the Flux subgroup $\Gamma = F(\pi_1(T^d))$ is discrete.
It turns out that every torus action can be reduced
to a cohomologically free torus action. More precisely,
we say that
a torus action
on $M$ is {\it reducible} if the following
conditions (i)-(iii) hold:

(i) there exists a $T^d$-invariant
symplectic
submanifold $N \subset M$ such that 
the incusion $H_1(N,\R) \to H_1(M,\R)$ is onto
(we use the convention that a point is a symplectic submanifold);

(ii) there exists a symplectic $T^r$-action on $N$ and
a surjective homomorphism $j: T^d \to T^r$ such that
$g|_N$ and $j(g)$ coincide as diffeomorphisms of $N$
for all $g \in T^d$;

(iii) $r+{\rm dim}\,N < d+ {\rm dim}\,M$.

This definition splits in two cases: either ${\rm dim}\,N < {\rm dim}\,M$
and in this case on can choose $j$ as the identity, or 
${\rm dim}\,N = {\rm dim}\,M$ which means that the whole action
factorizes through a smaller torus.

   As an example, the product of two standard circle actions
on $T^2 \times S^2$ is reducible since we may take 
$T^r = S^1 \times \1$ and $N$ to be
 the product of $T^2$ with a fixed point on $S^2$.

\begin{prop} \label{thm:CF} Every irreducible torus
action is cohomologically free.
\end{prop}
\proof{}
Assume that the action of $T^d$ on $M$ is not cohomologically
free. Then for some non-zero element $v$ in the Lie algebra
of $T^d$ the action generated by $v$ is Hamiltonian.
Denote by $H$ the Hamiltonian function.
Let $L \subset T^d$
be the subtorus (of positive rank!)
defined as the closure of the 1-parametric subgroup
$V$ generated by $v$. Notice that every critical point of $H$
is fixed by the action of $L$, since $V$ is dense in $L$.
Since the action of $L$ is linearizable near each fixed point,
we conclude that in some local complex coordinates $(z_1,...,z_n)$
on $M$ near a critical point one can write $H$ (up to an additive
constant) as
a linear combination of the $|z_i|^2$. In particular, either $H$ is
{\it identically constant} or
$H$ is {\it a Morse-Bott function with even indices}.

In the first case, the action of $L$ on $M$ is trivial.
Taking $N=M$, $T^r = T^d/L$, and defining $j$ as the natural
projection we find that the action is reducible.

In the second case,
take $N$ to be the minimum set of $H$. It follows
from the previous discussion that $N$ is a symplectic
submanifold, 
and  Morse-Bott theory implies that the
inlusion $H_1(N,\R) \to H_1(M,\R)$ is an epimorphism.
Moreover, since the action of $L$ fixes $N$ and commutes
with the action of the whole group $T^d$, we conclude
that $N$ is $T^d$-invariant. Therefore, taking $T^r = T^d$ and
defining $j$ as the identity map
we see that the action of $T^d$ is reducible.
This completes the proof.
\QED

\medskip

\NI
{\bf Proof of Theorem~\ref{thm:autonomous}.}

Consider a $T^d$- action on $M$.
Applying Proposition~\ref{thm:CF} repeatedly we end up with a
cohomologically free reduction since the process
terminates in view of (iii)!
In other words
there exists
a $T^d$-invariant
symplectic submanifold $N\subset M$ and
a homomorphism $j:T^d \to T^r$ which
satisfy conditions (i)-(iii) above and
such that the $T^r$- action on $N$
is cohomologically free.

Let $\Phi: H_1(T^r,\R) \to H^1(N,\R)$ and
$F:H_1(T^d,\R)\to H^1(M,\R)$ be the flux homomorphisms
of the $T^r$-action on $N$ and the $T^d$-action on $M$ respectively.
Let $j_*:H_1(T^d,\R) \to H_1(T^r,\R)$ and
$i^*:H^1(M,\R) \to H^1(N,\R)$ be the natural
maps. 
Clearly, $\Phi \circ j_* = i^* \circ F$, and thus
$$i^*(F(\pi_1(T^d))) = \Phi(j_*(\pi_1(T^d))).$$
The right hand side is discrete since
$j_*$ maps integer homology to integer homology
and $\Phi$ is a monomorphism as the flux of
a cohomologically free action. Moreover,
$i^*$ is a monomorphism in view of the condition (i).
Hence the flux subgroup
$\Gamma = F(\pi_1(T^d))$ is discrete.
This completes the proof.
\QED

\begin{cor} Every symplectic torus action on
a closed manifold splits into a product of a Hamiltonian
action and a cohomologically free action.
\end{cor}
\proof{} This
was proved by Ginzburg \cite{GI} for the case when either the
symplectic manifold is Lefschetz, or the symplectic form
represents an integer cohomology class. Combining
Ginzburg's argument with our
Theorem~\ref{thm:autonomous} one immediately gets
this statement for general manifolds.  \QED

\subsection{$S^1$ actions}

It was shown in~\cite{McD:mom} that a symplectic $S^1$-action on a
closed $4$-manifold $(M,\om)$ is Hamiltonian if and only if it has fixed
points.  Therefore, part (i) of Proposition~\ref{prop:orb} follows from
the next lemma.

\begin{lemma}  In the $4$-dimensional case every symplectic
$S^1$-action with contractible orbits has fixed points. 
\end{lemma} 
\proof{}
Consider
an action with no fixed points.  This is generated by a nonvanishing vector
field $X$.  As in~\cite{McD:mom}, one can slightly perturb $\om$ to an
invariant form which represents a rational cohomology class and then
rescale $\om$, to reduce to the case when the flux of the action, which is
represented by the form $i_X(\om)$, is integral.   In this case there is a
map $\mu: M\to S^1= \R/\Z$, called the generalised moment map, such that
$$
i_X(\om) = \mu^*(dt).
$$
Since $X$ has no fixed points $\mu$ is a fibration with fiber $Q$
and we may think of $M$ as made from $Q\times [0,1]$ by identifying
$Q\times \{0\}$ with $Q\times \{1\}$.  Moreover, the action of $S^1$ on
$Q\times [0,1]$ is Hamiltonian and so has the usual structure of such
actions.  Thus there is a Seifert fibration $\pi:
Q\to \Si$ and a family $\si_t$ of area forms on $\Si$ such that
$$
\om|_{Q\times \{t\}} = \pi^*(\si_t),\qquad [\si_t] - [\si_s] = -(t-s) c,\;\;s,t\in
[0,1], $$
 where $c\in H^2(\Si,\Q)$ is the Euler class of $\pi$.

We claim that it is impossible for the orbits to be contractible in
$M$.   To see this, observe that the long exact homotopy sequence for $Q\to
M\to S^1$ implies that $\pi_1(Q)  $ injects into $ \pi_1(M)$.  Hence if the
orbits are contractible in $M$, they  contract in $Q$.\footnote
{
This is the place where the argument fails in homology.  It is possible for
the orbit to represent a nonzero element of $H_1(Q)$ while being
nullhomologous in $M$.}
    But if they contract
in $Q$, the Euler class $c$  must be
nontrivial.  Therefore $[\si_0]\ne [\si_1]$, and it is impossible to find a
suitable gluing map $Q\times \{0\}\to Q\times \{1\}$ since such a map
would induce a symplectomorphism $(\Si,\si_0) \to (\Si,\si_1)$.\QED

\section{Proof of the main Theorem in the $C^1$-case} \label{se:C1}

We will begin by proving the elementary results stated in
Propositions~\ref{prop:equivalence} and~\ref{prop:commutative}. 
We then prove the various parts  of  Theorem~\ref{thm:main-C1},
modulo a lemma needed for (ii) that is presented in \S3.2.

\subsection{The main arguments}

\NI
{\bf Proof of Proposition~\ref {prop:equivalence}.} 

We must show that the $C^1$ closure of $\Ham(M)$ is equivalent to the
discreteness of the flux group $\Ga$.  

Let $\phi_k$ be a sequence of
Hamiltonian diffeomorphisms which $C^1$-converges to some symplectic
diffeomorphism $\phi$. By composing  the sequence 
and its limit with $\phi_{k_0}^{-1}$ for some fixed $k_0$,
we get a new sequence $\psi_k$ of Hamiltonian diffeomorphisms converging
to some $\psi$ which is Hamiltonian iff $\phi$ is. By choosing $k_0$
large enough, the new sequence 
is $C^1$-close to the identity for $k > k_0$, 
and therefore
there are closed $1$-forms $\la, \la_k$ on $M$ such 
that the
autonomous symplectic path $\psi^{\la}_{t \in [0,1]}, 
\psi^{\la_k}_{t \in [0,1]}$ generated 
by the $\om$-duals of $\la, \la_k$ are small paths whose endpoints are 
$\psi, \psi_k$
and whose fluxes are equal to $[\la], [\la_k]$. For each $k\ge k_0$,
we thus get a loop by composing (in the sense of 
paths, not pointwise) a  Hamiltonian
path $\{\phi^k_t\}$ from the identity to $\phi^k_1 = \psi^k$ with the path
$\{\psi^{\la_k}_{1-t}\}$. This loop has flux $-[\la_k]$, which is 
arbitrarily small
for  $k\ge k_0$.  If $\Ga$ is discrete, these classes must vanish,
and therefore their limit $[\la] = 0$. But this means that the path
$\psi^{\la}_{t \in [0,1]}$ is Hamiltonian.  Hence $\psi^\la_1 = \psi$
is Hamiltonian, as required.

   Conversely, if $\Ga$ is not discrete, 
there is an arbitrarily small class
$[\la] \in H^1(M;\R)$ which is not in $\Ga$ but is the limit of classes
in $\Ga$. Then the above construction 
 gives a sequence of
Hamiltonian diffeomorphisms 
which $C^1$-converge to a nonHamiltonian one.
\QED
\MS

\NI
{\bf Proof of Proposition \ref {prop:commutative}.}

We must show that, if the flux homomorphism $F$ has
discrete image on the subgroup formed by
symplectic loops with  trivial
evaluation in $\pi_1(M)$, its whole image $\Ga$ is discrete.  
To do this, it is clearly enough
to show that there exists a neighborhood $U$ of $0$
in $H^1(M,\R)$ with the following property:
\begin{quote}
\em
 every class $a \in \Gamma \cap U$ is the image under  $F$ of
 a symplectic loop with contractible orbits.
\end{quote}

If $U$ is sufficiently small, then every class $a\in U$ can be
represented by a
$1$-form $\la$ which generates an autonomous locally-Hamiltonian 
flow $\psi_t^\la$
whose only closed orbits in the time interval $t \in
[0;1]$ are the zeros of $\la$.  
We claim that such $U$ have the required property.

To see this
suppose that $a \in \Gamma$. Then $\psi_1^\la$
is the time-$1$ map of some Hamiltonian flow 
$\{\phi_t\}$.  By the discussion
 after Proposition~\ref{prop:commutative} there exists a fixed
point $x$ of $\phi_1 = \psi_1^\la$ such that the corresponding
closed orbit of $\{\phi_t\}$ is contractible.
As in the proof of Proposition~\ref{prop:equivalence}, now
consider  the loop in $\Symp_0(M)$ obtained by composing the paths
$\{\phi_{1-t}\}$ and $\{\psi^{\lambda}_t\}$. We see that its flux equals $a$.  
Moreover, if
the only fixed points of $\psi_1^\la$ are zeros of $\la$, its orbits are 
contractible.
This completes the proof.
\QED
\MS

\NI
{\bf Proof of Theorem \ref{thm:main-C1} (i)}

Assume that the manifold $M$ is spherically rational.
We must show that $\Ga$ is discrete. 
In view of Proposition~\ref{prop:commutative} it suffices
to consider symplectic loops with contractible orbits only.
Let $\{\phi_t\}$ be such a loop, and let $a \in H^1(M,\R)$
be its flux. Then the value of $a$  
on each closed curve $\gamma \subset M$
is equal to the symplectic area of the torus $C_\ga$.
But this torus is the image of a sphere
since one of its generating loops is contractible in $M$. Thus the 
homomorphism 
$$
H_1(M, \R) \stackrel{ev_\phi}{ \longrightarrow} H_2(M, \R) \stackrel{\int
\om}{ \longrightarrow}  \R
$$
 factors through spheres 
$$
a:H_1(M, \R)  \stackrel{ev_\phi}{ \longrightarrow} H^S_2(M, \R) 
\stackrel{\int \om}{ \longrightarrow} \R.
$$
 But if the manifold is spherically
rational, this map takes values in some discrete subgroup
of $\R$ which does not depend on a particular choice of
$a$. This implies the discreteness of the flux subgroup $\Gamma$.
\QED
\MS

\NI
{\bf Proof of Theorem \ref{thm:main-C1} (ii)}

We must show  that if the minimal Chern number of $(M, \om)$ is
zero or is $\ge 2n$, then the Flux conjecture holds. This is
based on the comparison of the Maslov and Floer indices. 

In view of Proposition~\ref{prop:commutative} 
it is enough to  restrict
our attention to the subgroup $K$ of $\pi_1(G)$
consisting of  loops 
whose evaluation in $\pi_1$ vanishes. 
Fix a Riemannian metric $g$ on $M$, and define a norm $||a||$
of a class $a \in H^1(M,\R)$ as the $C^1$-$g$-norm
of the unique $g$-harmonic form in this class.
The following lemma is proved in \S~3.2. 

\begin{lemma} \label{le:1}  There exists $\epsilon > 0$ such that every
nonzero cohomology class $a$ with $||a|| < \epsilon$
is represented by a Morse form $\la$ which has no critical
point of index $0$ and $2n$, and whose Hamiltonian
flow has no nonconstant $1$\/-periodic orbit.
\end{lemma}

Granted this, assume that the image of $F(K)$ is not 
discrete, and choose $\epsilon$ as in Lemma~\ref{le:1}.
Then there exists $\be \in K$ with
$0 < ||F(\be)|| < \epsilon$.
Represent the class $F(\be)$ by a form $\la$
as in Lemma~\ref{le:1} and denote by $\psi_1 $ the time-$1$
map of  the corresponding
symplectic flow $\{\psi_t^\la\}$. In view of our choice, the only $1$-periodic
orbits of $\{\psi_t^\la\}$ are the constant ones.

 Since $\psi_1$ is a Hamiltonian diffeomorphism,
we can join it to the identity by a Hamiltonian path $\{\phi_t\}$.
Introduce the following notations. For a fixed point $x$ of $\phi_1=\psi_1$
denote by $i_\la(x)$ the Conley--Zehnder index coming from the flow of
$\la$ and by $i_H(x)$ the Conley--Zehnder
 index coming from the Hamiltonian path.
If the minimal spherical Chern number $N$ is $0$ or
is $\ge n$, then  Hofer-Salamon proved in~\cite{HS} that
the Floer homology of a Hamiltonian path is isomorphic to the 
ordinary homology
with coefficients in the Novikov ring of the group $H_2^S(M,\Z)$. In
particular there are fixed points $x,y$ of $\phi$ with
$$
i_H(x) - i_H(y) = 2n \quad \pmod {2N}
$$
(Recall that all relative indices of Floer homology are
defined only modulo $2N$).
On the other hand the orbits of the loop which is the composition of 
$\{\phi_t\}$
with $\{\psi_{1-t}^\la\}$ have indices which are independent of the 
choice of the orbit, denote them by the constant
$m \in \Z$ modulo $2N$.
Thus, for any fixed point $z \in M$,
$$ 
i_H(z) - i_\la(z) = m \quad \pmod {2N}
$$ 
and therefore
$$
i_H(x) -i_H(y) =  i_\la(x) - i_\la(y) \quad \pmod {2N},
$$
which implies that
$$
i_\la(x) - i_\la(y) =  2n \quad \pmod {2N}.
$$
But  our choice of $\la$ implies that 
$$
-2n+2 \le i_\la(x) - i_\la(y)\le 2n-2
$$ 
which means that $2\le 2n -(i_\la(x)-i_\la(y))\le 4n-2$. The last equation states
that $2n -(i_\la(x)-i_\la(y))$ is a multiple of $2N$. But this is impossible
if $N$ is $0$ or $\ge 2n$.
\QED

\NI
{\bf Proof of Theorem~\ref{thm:main-C1}(iii).}

Suppose, by contradiction, that the flux
is not discrete when $\pi_1$ acts trivially on $\pi_2$ and $M$ has
dimension $n=4$. Because it is non-discrete, it is not factorizable. Thus
there is a loop with trivial evaluation in $\pi_1(M)$ and non-zero
flux $a \in H^1(M)$. By the commutativity of the diagram
in \S~1.1~(i), we must have $a \cup [\om^{n-1}]=0$. Hence there is
a non zero class such that:
\begin{description}
\item
(i)
$a\cup[\om] = 0$; and 

\item
(ii) there is a loop of symplectic
diffeomorphisms $\{\phi_{t \in [0,1]}\}$ whose flux is $a$ and which
has contractible evaluation in $\pi_1(M)$. 
\end{description}

Now choose an element $\ga \in \pi_1M$
such that $a(\ga) = \ka \neq 0$, and let $C =C_\ga$, the trace of $\{\phi_t\}$
on $\ga$.
Then, if $p: \Tilde M\to M$ is the universal cover, 
$C$ lifts to a $2$-sphere $\Tilde C$ in $\Tilde M$, and
$$
p^*([\om])(\Tilde C) = [\om](C) = \ka \ne 0.
$$
 Therefore it suffices to prove the following lemma.\footnote
{
This lemma gives some results in all dimensions.  However, it is only in
dimension $4$ that one gets a clean statement about the Flux conjecture.}

 \begin{lemma}   Let $a$ be a nonzero class in $H^1(M;\R)$ such
that for some $k$ 
$$
a\cup [\om]^k = 0,\qquad p^*([\om]^k) \ne 0.
$$
Then $\pi_1$ acts nontrivially on $H_{2k}(\tilde M)$.  Moreover, if
$p^*([\om]^k)$ does not vanish on some element of
$H_{2k}^S(\tilde M)$, $\pi_1(M)$ acts nontrivially on $\pi_{2k}(M)$.
\end{lemma}

\proof{}
Choose an element $\ga \in \pi_1M$
such that $a(\ga) = \ka \ne 0$, and let 
 $\tau$ be the deck
transformation associated to $\ga$.  Then, if the $1$-form $\la$ represents
$a$, we have $p^{\ast}(\la)= dH$ where $\tau^{\ast}H - H = \ka$.  Further,
$$ 
\la \wedge \om^k = d \theta
$$
for some $2k$\/-form $\theta$ on $M$ because $a\cup [\om]^k = 0$.
This implies
$$
d(H p^{\ast} \om^k) = d(p^{\ast}\theta),
$$
which implies  in turn that
$$
H p^{\ast} \om^k = p^{\ast} \theta + \si
$$
for some closed $2k$\/-form $\si$ on $\Tilde M$. Hence:
$$
\tau^{\ast}\si - \si = (\tau^{\ast} H - H)p^{\ast} \om^k = \ka p^{\ast}\om^k,
$$
and so there is a $2k$-cycle $\Tilde C$ in $\Tilde M$ such that
$$
\int_{\tau(\tilde{C})} \si - \int_{\tilde{C}} \si = 
\int_{\tilde{C}} (\tau^{\ast}\si - \si) = \ka \int_{\Tilde C} p^*\om^k 
\neq 0. $$
Therefore $\tau(\Tilde{C})$ and $\Tilde{C}$ are not homologous in
$\Tilde M$, and the action
of $\pi_1(M)$ on $H_{2k}(\tilde M)$ is nontrivial. 
Finally, if $p^{\ast}([\om]^k)$
does not vanish on some element of $H_{2k}^S(M)$,  the same argument
shows that the action of $\pi_1(M)$
on $\pi_{2k}(\tilde M)$ is nontrivial too (note that this action given
by the deck transformations is well-defined because $\pi_1(\tilde M)$
vanishes). Because $\pi_{2k}(\tilde M)$ is isomorphic to $\pi_{2k}(M)$,
we conclude that $\pi_1M$ acts nontrivially on $\pi_{2k}(M)$.
\QED

\subsection{Proof of Lemma~3.1} \label{ss:harmonic}

It was shown 
by G. Levitt \cite{Le}  
that every non-zero cohomology class can be
represented by a Morse 1-form
without critical points of index $0$ and $2n$.
We show here that the lemma~3.1 follows easily from
this statement. However one can also use a different approach
based on the fact that a generic harmonic
1-form is Morse (see the Appendix).

Here is a short proof of Lemma~3.1 based on Levitt's theorem. 
On some small closed convex neighbourhood $V$ of
$0 \in H^1(M,\R)$ choose any smooth section $\si$ of $Z^1(M) \to H^1(M, \R)$
(where $Z^1(M)$ is the space of closed 1-forms) which  is $0$ at $0$
(one can do this by choosing a Riemannian metric say).
Then take any smooth convex  hypersphere $S$ in $H^1(M,\R)$
containing $0$ in its interior, and for each class $[\la]$  on $S$,
choose a Levitt form $\la$ (which means here a Morse form with no 
critical point of index $0$ or $2n$) and
define the neighbourhood $U_{\la} \subset S$ as the intersection of
$S$ with the classes $[\la] + rV$, where $r$ is a small real number. Each
element in $U_{\la}$ is represented by the form $\la + r \si(v)$. Thus
if $r$ is small enough, all elements of $U_{\la}$ are represented by
Levitt forms and this representation is continuous. Now there is a finite
covering of $S$ by such subsets. This implies that there is $\eps$
small enough so that the hypersphere $\eps S$ as well as its interior
is covered by Levitt forms
that are such that the Hamiltonian flow has no non-constant 1-periodic orbit.

\section{The $C^0$-conjecture and Lefschetz 
manifolds} \label{se:C0}

In this section we prove Theorem \ref{thm:Lefschetz}.
We start with a discussion of the properties of
the flux homomorphism in the more general context
of paths of {\it maps}. 

Let $M$ be a closed manifold, and let
$\al= \{\phi_t\}_{t\in [0;1]}$ be a  path of smooth maps $M\to M$
such that $\phi_0 = \1$ and $\phi_1$ is a diffeomorphism. Let
$I_\al^k$ be the ring of $\phi_1$-invariant differential forms on $M$
of degree $k$. 
Define the flux map $\Phi : I_\al^k \to H^{k-1}(M,\R)$
as follows. Given an invariant form $\sigma \in I_\al^k$ and
a $(k-1)$-cycle $\gamma$ on $M$, the value $\Phi(\sigma)$
on $\gamma$ is defined to be the integral of $\sigma$ on the trace
$\cup_{t\in [0,1]}\,\{\phi_t(\ga)\}$ of $\gamma$ under the path. It is not hard to
see that this value depends only on the homology class of $\gamma$,
so that $\Phi$ is well defined.  Note that its domain $I_\al^k$ depends on
$\phi_1$ in general, and may well be zero when $\phi_1\ne \1$.

It turns out that $\Phi$ inherits the important {\it derivation}
property of the usual flux associated to a closed path
of maps which starts and ends at the identity. Namely,
given two $\phi_1$-invariant
forms $\sigma_1,\sigma_2$ (which for convenience
we assume to  have  even degrees) we have
$$
\Phi(\sigma_1 \wedge \sigma_2) = [\sigma_1] \cup \Phi (\sigma_2) +
\Phi[\sigma_1] \cup [\sigma_2].
$$
To see this, suppose first that the maps $\phi_t$ are diffeomorphisms,
generated by the family of vector fields $X_t= \dot{\phi}_t$.  
Then it is easy to see that
 the class $\Phi(\si)$ is represented by the form 
$$
\int_0^1\phi_t^*(i_{X_t}\si)\,dt = \int_0^1i_{Y_t}(\phi_t^*\,\si)\,dt,
$$
where $Y_t = \phi_t^*(X_t)$.  With this representation, the derivation
property is obvious. In the general case it is not hard to make sense of the
right hand side of this formula, because $Y_t$ can be considered as a vector
field along $\phi_t$ (when $\phi_t$ is not a diffeomorphism).

In what follows we write $ \Phi =\Phi_{\alpha}$ in order
to emphasize the dependence of $\Phi$ on the path $\alpha$.

Return now to our symplectic situation. Let $(M,\omega)$
be a closed symplectic manifold. Denote by $\Xx$
the set of paths of maps $\alpha= \{\phi_{t \in [0;1]}\}$ 
with symplectic endpoints:
$\phi_0 = \1$
and $\phi_1 \in \Symp_0(M)$. We will endow $\Xx$ with the
$C^0$-topology associated to a Riemannian metric on $M$.

The proof of Theorem \ref{thm:Lefschetz} is based on the following
statement.

\begin{lemma} \label{le:CONT} Assume that $M$ is Lefschetz.
Then the map 
$$\Xx \to H^1(M,\R),\;\; \alpha \to \Phi_{\alpha}(\omega)$$
is continuous.
\end{lemma}
\proof{}The derivation property implies that
$\Phi_{\alpha}(\omega^n) = n[\omega]^{n-1} \cup \Phi_{\alpha}(\omega)$.
The mapping $\alpha \to \Phi_{\alpha}(\omega^n)$ is continuous
since the volume of a top-dimensional subset of $M$
varies continuously under continuous deformations of the subset.
The needed statement follows now from the Lefschetz property.
\QED

\MS

\NI
{\bf Proof of Theorem \ref{thm:Lefschetz}.}

Suppose that the sequence $\psi_n \in \Ham(M)$ converges
$C^0$ to a symplectomorphism $\psi$ in the identity
component $\Symp_0(M)$. We have to show that the flux $F(\psi)$
belongs to the quotient group $\Gamma_{top}/\Gamma$. 

Fix $n$ large enough
and 
consider the map
$\phi = \psi_n^{-1}\psi$. Define two paths of maps between
the identity and $\phi$. One of them, say $\alpha$, is a short
path of maps along
geodesics in $M$. (If $n$ is large then $x$ and $\phi( x)$ can be
joined by a unique geodesic in $M$ for every point  $x \in M$.)
The second path, say $\alpha'$, is a path of symplectic 
diffeomorphisms. Let $\beta$ be the loop of maps formed by
composing of $\alpha'$ with $-\alpha$, where $(-\al)_t = \al_{1-t}$. Clearly,
$$
F(\beta) = F(\alpha')-\Phi_{\alpha}(\omega).
$$ 

Note that $F(\beta)$ belongs to $\Gamma_{top}$,
and that $F(\alpha') = F(\psi) \;\;({\rm mod}\;\Gamma)$. Further,
Lemma \ref{le:CONT} implies that
$\Phi_{\alpha}(\omega)$ can be made arbitrarily small
by choosing $n$ sufficiently large. But both $\Gamma$ and
$\Gamma_{top}$ are discrete since $M$ is Lefschetz.
Therefore we get that $F(\psi) = F(\beta)\;\;({\rm mod}\;\Gamma)$.
This completes the proof.
\QED

\section{Appendix}

We present here another more geometric proof of Lemma~\ref{le:1}, based on 
the representation of forms by harmonic forms. We first need the 
following folkloric 
statement, that was mentioned to us
by M. Farber and M. Braverman.  

\begin{lemma}\label{le:2} After an arbitrarily small smooth
perturbation of the metric, the harmonic form in a given
nonzero cohomology class is Morse.
\end{lemma}

We prove this by a method adapted from Uhlenbeck \cite{UHL}.

\SS
\proof{}
Let $\cal R$ be the space of all Riemannian metrics on $M$,
and let $\Omega$ be the space of all closed 1-forms on $M$
representing a given nontrivial cohomology class in $H^1(M,\R)$.
Define the subset $\Lambda \subset {\cal R} \times \Omega$ consisting
of all pairs $(g,\lambda)$ such that $\lambda$ is $g$-harmonic.

The sets $\cal R$, $\Omega$ and $\Lambda$ have a natural structure
of Banach manifolds with respect to suitable Sobolev norms.
Moreover, the natural projection $\Lambda \to \cal R$ is a smooth
one-to-one map.
Hence in view of Smale-Sard transversality theorem with parameters,
it is enough to check that the evaluation map
$A : \Lambda \times M \to T^*M$ defined by
$$A(g,\lambda,x) = (x,\lambda_x)$$
is {\it transversal to the zero section}.

Assume that $\lambda$ is a $g$-harmonic form which vanishes
at some point $q$. Note that $T_{(q,0)}T^*M = T^*_qM \oplus T_qM$.
In order to prove the needed transversality result it suffices to
show the following:
\begin{quote}{\em
for every $\xi \in T^*_qM$ there is a tangent
vector 
$$
(\dot g, \dot \lambda) \in T_{(g,\lambda)}\Lambda
$$
such that $\dot \lambda_q = \xi$.
}
\end{quote}

First of all we  compute the tangent space
$T_{(g,\lambda)}\Lambda$ that we consider as a linear subspace
in the product $S^2(M) \times \Omega_0(M)$,
where $S^2(M)$ is the space of all symmetric 2-forms on $M$
and $\Omega_0(M)$ is the space of
all exact 1-forms.

Denote by $g^{-1}$ the natural operator associated
with the metric $g$ which
transforms $(0,k)$-tensors to $(1,k-1)$-tensors.
A straightforward computation shows that a
vector $(\dot g, \dot \lambda) \in S^2(M) \times \Omega_0(M)$
is tangent to $\Lambda$ at $(g,\lambda)$
if and only if:
$$
\Delta_g \dot \lambda - d\Phi(\dot g) = 0,
$$
where
$$
\Phi(\dot g) = {\rm div}_g (-(g^{-1}\dot g) g^{-1}\lambda
+ {1 \over 2} {\rm tr}(g^{-1}\dot g) g^{-1}\lambda ).
$$

Since $\la$ does not vanish identically, we can find an
open set $B \subset M$ such that $\lambda$ has
no zero in $B$.
For every smooth function $f$
which is supported in $B$ and has zero $g$-mean,
one can find a symmetric form (also supported in $B$)
$\dot g$ in such a way that $\Phi(\dot g) = f$.
Also,
$\Delta_g \dot \lambda = d \Delta_g u$ where $du = \dot \lambda$,
and hence $(\dot g, \dot \lambda)$ is tangent to $\Lambda$
at $(g, \lambda)$ if $\Delta_g u = f $.
In view of this discussion, the required
transversality fact follows immediately from the next lemma.
We write $<f >$ for the $g$-mean value of a function
$f$.

\medskip

\begin{lemma} Let $(M,g)$ be a closed Riemannian
manifold, $B \subset M$ be an open subset and $q \in M$
be a point which is chosen away from the closure of $B$.
Then for each $\xi \in T^*_qM$ there exists a smooth
function $f$ supported in $B$ with $<f > = 0$
such that the solution of the equation $\Delta_g u = f$
satisfies $d_qu = \xi$.
\end{lemma}

\proof{}
Define the operator $E : C^{\infty}(M) \to C^{\infty}(M)$
such that $<Ef > = 0$ and
$$
\Delta_g Ef = E\Delta _g f = f -<f >.
$$
Its kernel (or {\it Green's function}) $e(x,y)$ is defined by
$$(E f )(x) = \int_M e(x,y) f(y) dy \; ,$$ and satisfies the
equation $\Delta_g e(x,.) = -1$ on $M - \{ x \}$.

Introduce local coordinates $x_1, ... ,x_n$ near $q$ and
identify tangent and cotangent spaces to $M$ at $q$ with $\R^n$.

Let $V$ be the space of all smooth functions with zero mean
which are supported in $B$, and define an operator
$T : V \to \R^n$ by
$$
T(f ) = {\partial (Ef ) \over \partial x} (q).
$$
We have to prove that $T$ is surjective. Assume, on the contrary,
that there exists a vector $w \in \R^n$ such that $(w, Tf ) = 0$
for all $f \in V$ (here  $(\cdot, \cdot)$ is a scalar product).
Then we immediately get that the function $b : M -\{ q\} \to \R$,
$$
b(y) = (w, { \partial e \over \partial x } (q,y)), 
$$
is {\it constant} on $B$. Notice that $b$ is harmonic, and
in view of elliptic regularity it is therefore constant everywhere
on $M -\{q\}$.

Let us show that this is impossible. Indeed choose
a function $h$ with $<h > = 0$ which is harmonic near $q$
and such that $(w, {\partial h \over \partial x} (q)) \neq 0$.
We write 
$$
h (x) = \int_M e(x,y) \Delta_g h (y) dy.
$$
Since $\Delta_g h (y)$ vanishes near $q$, the derivative
of the integral at $q$ equals the integral of the derivative,
and hence
$$
(w, {\partial h \over \partial x} (q)) =
\int_M b(y) \Delta_g h (y) dy = 0
$$
since $b$ is constant. This contradiction proves the lemma.
\QED

\NI
{\bf Proof of Lemma \ref{le:1}} 

Fix $\epsilon > 0$ such that the Hamiltonian
flow of every closed 1-form $\la$ with $C^1$-$g$-norm $< 2\epsilon$
has no nonconstant 1-periodic orbits. Fix any nontrivial class
of norm $< \epsilon$.
Take a very small perturbation 
$g'$ of the metric $g$ such that the $g'$\/-harmonic form
$\la$ in this class
is Morse and has $C^1$-$g$\/-norm less than $2\epsilon$.
Notice that a harmonic Morse form has no critical points of index
0 and $2n$ by the maximum principle.
This completes the proof.
\QED


\begin{thebibliography}{9999}

\bibitem{B} 
A. Banyaga,   Sur la structure du groupe des
diff\'eo\-mor\-phismes qui pr\'e\-ser\-vent une forme symplectique.
{\it Commentarii Mathematicae Helveticae\/}, {\bf 53},
(1978) 174--227. 

\bibitem{BPGAFA}
M. Bialy and L. Polterovich, Hamiltonian diffeomorphisms
and Lagrangian distributions, {\it Geom. and Funct. Analysis},
{\bf 2} (1992), 173-210.

\bibitem{CZ} 
C. Conley and E. Zehnder,  
The Birkhoff-Lewis fixed point theorem and a conjecture of V.I. Arnold,
{\it Inventiones
Mathematicae\/}, {\bf 73}, (1983) 33--49. 


\bibitem{FO}  K. Fukaya and K. Ono, Arnold conjecture and
Gromov--Witten invariants, preprint (1996)

\bibitem{GI}  V. L. Ginzburg, Some remarks on symplectic 
actions of compact groups, {\it Math. Zeitschrift}, {\bf 210} (1992) 625--640.

\bibitem{HS} 
H. Hofer and D.A. Salamon,  Floer homology and Novikov rings,  
{\it Floer Memorial volume\/}  eds. Hofer,
Taubes, Weinstein, Zehnder, Birkh\"auser, Basel, (1996).

\bibitem{LMP} F. Lalonde, D. McDuff and L. Polterovich, The topological
rigidity of Hamiltonian loops, in preparation.

\bibitem{Le} G. Levitt, 1-formes ferm\'ees singuli\`eres et groupe
fondamental, {\it Invent. Math.}, {\bf 88} (1987), 635-667.

\bibitem{LIUT} G. Liu and G. Tian, Floer homology and Arnold
conjecture, preprint (1996)


\bibitem{LO}
G. Lupton and J. Oprea,  Cohomologically symplectic spaces, 
Toral actions and the Gottlieb group, {\it Trans. AMS} {\bf 347} (1995),
261-288.

\bibitem{McDF} 
D. McDuff,   Symplectic diffeomorphisms 
and the flux homomorphism. {\it Inventiones Mathematicae}, {\bf 77},
(1984) 353--66. 

\bibitem{McD:mom} 
D. McDuff,  The moment map for circle actions on symplectic manifolds,
{\it Journal of Geometrical Physics\/}, {\bf 5} (1988) 149--60. 

\bibitem{MS}  D. McDuff and D.A. Salamon, {\it Introduction to
Symplectic Topology},  OUP, Oxford, (1995)

\bibitem{SEI} P. Seidel, $\pi_1$  of symplectic automorphism groups
and invertibles in quantum cohomology rings, preprint dg-ga/9511011,
(1995)

\bibitem{TH}  D. Th\'eret, in preparation.

\bibitem{UHL}  K. Uhlenbeck,
Generic properties of 
eigenfunctions, {\it Amer. J. of Mathematics}, {\bf 98} (1976), 
1059-1078.


\end{thebibliography}
\end{document}